\documentclass[12pt,reqno]{amsart}
\usepackage{amsmath,amsxtra,latexsym,amsthm,amssymb,amscd,pb-diagram}
\usepackage{mathrsfs,mathabx,amsfonts}
\usepackage[colorlinks=true,linkcolor=red,citecolor=blue]{hyperref}
\usepackage[margin=1in]{geometry}
\usepackage{bm}
\usepackage{color}
\usepackage{makecell,multirow,diagbox}
\usepackage{mathtools}
\usepackage[displaymath, mathlines]{lineno}
\usepackage{dsfont}
\usepackage{graphicx}
\usepackage{epsfig}
\usepackage{array}
\usepackage{enumitem}

\newtheorem{theorem}{Theorem}[section]
\newtheorem{lemma}[theorem]{Lemma}

\newtheorem{corollary}[theorem]{Corollary}

\title[Wave system on weighted graphs]{Non-existence results for a system of wave inequalities on locally finite graphs}

\author[A. T. Duong]{Anh Tuan Duong}
\address{Anh Tuan Duong\\ Faculty of  Mathematics and Informatics, Hanoi University of Science and Technology, 1 Dai Co Viet, Bach Mai, Ha Noi, Viet Nam.}
\email{tuan.duonganh@hust.edu.vn}

\author[T.A. Dao]{Tuan Anh Dao}
\address{Tuan Anh Dao\\ Faculty of  Mathematics and Informatics, Hanoi University of Science and Technology, 1 Dai Co Viet, Bach Mai, Ha Noi, Viet Nam.}
\email{dao.tuananh@hust.edu.vn}

\usepackage{color}

\begin{document}
\begin{abstract}
	Let $V$ be a locally finite, connected and weighted graph. We study   non-existence results of non-trivial, non-negative solutions of the system 
$$	\begin{cases}
		u_{t t}-\Delta u \geq  h_1|v|^p & \text { in } V \times(0, \infty), \\ 
		v_{t t}-\Delta v \geq h_2|u|^q& \text { in } V \times(0, \infty), \\ 
		u=u_0;\;v=v_0 & \text { in } V \times\{0\}, \\
		u_t=u_1;\;v_t=v_1  & \text { in } V \times\{0\},
	\end{cases}$$
	where $p,q>1$, $h_1,h_2$ are positive potentials. Under some volume growth condition of a ball, we prove that the system has no non-trivial non-negative solutions. In particular, our result is a natural extension of that in [\textit{D.~D.~Monticelli, F.~Punzo, and J.~Somaglia. Nonexistence results for the semilinear wave equation on graphs. arXiv.2506.08697, 2025.}] from a single inequality to a system.
\end{abstract}
\textbf
\maketitle

\section{Introduction}

In recent years, partial differential equations  on  weighted graphs have been extensively studied, see e.g \cite{Gr18,Ge18, GHJ18,GLY16,GLY16b,HL21,LW17,LW18,LY20,GHS23, MW13,IBR23,LZ19,HS21,HSZ20, Liu23,Liu24,MPS23,MPS24,BMP23,Wu21,Wu24,MDQ25,XW23,MPS26}. One of the most interesting subjects is Liouville-type theorem, i.e the non-existence of non-trivial solutions of partial differential equations \cite{GHS23, MPS23,MPS24, XW23, LW17,LW18} and references given there.

 On the class of elliptic equations,  the Lane-Emden inequality on weighted graphs 
 \begin{equation}\label{eat1}
 	\Delta u\geq u^p,
 \end{equation}
 was studied by Gu, Huang and Sun \cite{GHS23}. In \cite{GHS23}, the optimal Liouville-type theorem for positive solutions of \eqref{eat1} on weighted graphs was established under some $(p_0)$ condition and the growth of volume of a ball. Later, by removing the $(p_0)$ condition and adding the  existence of a pseudo-metric, the authors in \cite{MPS23} also proved the existence and non-existence of positive solutions of  \eqref{eat1}. In \cite{XW23}, the authors generalized the result in \cite{GHS23} to the case of $p$-Laplace operator. Recently, the results in \cite{GHS23,MPS23} have been generalized to a system of elliptic inequalities in \cite{MD23,MDN24}.

On the class of  parabolic equations on weighted graphs, we would like to mention the articles  \cite{LW17,LW18, Liu23, Liu24, Wu21,Wu24,DF24, MPS24} and references given there. The existence of global solutions of parabolic problems were investigated in \cite{LW17,LW18,Wu21,Wu24}. 
In  \cite{MPS24}, the authors have studied  the parabolic inequality on weighted graphs
 \begin{equation}\label{eat2}
	u_t-\Delta u\geq v(x) u^p,
\end{equation}
 with $p>1$  and the  potential $v$.  The non-existence of positive solutions of  \eqref{eat2} has been proved under the assumption of existence of a pseudo-metric. In addition, the sharpness is also discussed.  Remark that it is  a challenging problem in  studying the non-existence of positive solutions of \eqref{eat2} under the $(p_0)$ condition because of the presence of the time variable. Indeed, for \eqref{eat2},  the $(p_0)$condition does not imply any inequality of Harnack type which is a crucial tool in \cite{GHS23,MD23}. Very recently, some non-existence results for system of parabolic inequalities have been obtained in \cite{DF25}.

In contrast to the parabolic or elliptic equations, there have only a few works dealing with the wave equations on weighted graphs. Let us first mention some articles studying wave equations/system in Euclidean spaces or in Riemannian manifolds \cite{Ka80, PV00,  Zhang01, TY01, MP01b, JMY20, JSV21,  MPS20, JRST25}. Concerning the  wave inequality in Euclidean spaces,
\begin{equation}\label{ew1}
	u_{t t}-\Delta u \geq|u|^p,
\end{equation}
it was proved in \cite{Ka80}  that if the initial data satisfy some suitable positivity conditions, are compactly supported, and
\begin{equation}\label{ew8}
	1<p \leq \frac{N+1}{N-1}
\end{equation}
then \eqref{ew1} has no solution on $\mathbb{R}^N \times(0, \infty)$.  Later, this exponent was shown to be optimal in the paper \cite{PV00}. 

Concerning the equation $$
u_{t t}-\Delta u =|u|^p,
$$ the  non-existence of global solutions was studied in \cite{Si84} with the condition of exponent
$1<p<p^*,$
where $p^*$ is the positive root of the algebraic equation
$$
(N-1) x^2-(N+1) x-2=0.
$$
For the wave equation with  a linear damping term , i.e 
$$
u_{t t}+u_t-\Delta u =|u|^p,
$$ 
the authors in  \cite{TY01, Zhang01} proved that the influence of the damping is powerful enough to shift the critical exponent $p^*$ to the left. Moreover, the critical exponent in this case is the same as that in the parabolic one: $p_c(N)=\frac{N+2}{N}$.

The study of blow-up for wave equation in exterior domains 
$$
u_{t t}-\Delta u =|u|^p
$$ 
 was initiated  in \cite{Zhang01b}. Recently, the non-existence of global solutions in exterior domains of system of inhomogeneous wave inequalities 
 $$\begin{cases}
 	u_{t t}-\Delta u \geq|x|^a|v|^p,\\
 		v_{t t}-\Delta v \geq|x|^b|v|^q,
 \end{cases}$$
 has been obtained in \cite{JMY20}.
 
 In graph setting,  the non-existence of global solutions of the wave inequality  with potential
 \begin{equation}\label{eat3}
\begin{cases}
	 	u_{t t}-\Delta u \geq h(x)|u|^p, t>0,\\
	 	u(x,0)=u_0(x); u_t(x,0)=u_1(x),
\end{cases}
 \end{equation}
 has been recently studied in \cite{MPS26}, where the non-existence of non-trivial non-negative solutions was established under some volume growth conditions. Inspired by \cite{MPS26}, in this article, we propose to study the non-existence of non-trivial non-negative solutions of a system of wave inequalities involving potentials.
 
The paper is organized as follows. In the next section, we recall  the graph setting and present main results. The proof is given in the last sections.
\section{Graph setting and Main results}\label{s0}
Let $V$ be a countably infinite set whose elements are called vertices.  Define the edge weight $$ \omega: V\times V\to [0,\infty)$$ satisfying the following properties:

$\bullet $ (Absence of loops) $\omega_{xx}=0$ for all $x \in V$;

$\bullet $  (Symmetry of the graph) $\omega_{xy}=\omega_{yx}$ for all $(x, y) \in V \times V$;

$\bullet $ (Finite sum) $\sum\limits_{y \in V} \omega_{xy}<\infty$ for all $x \in V$.

\noindent The vertex weight is defined as a positive function $\mu: V\to(0,\infty)$. We shall present the relation between $\omega$ and $\mu$ later on.

\noindent Let  $x\in V$, we write  $y\sim x$ if $\omega_{xy}>0$ and in this case $x$ and  $y$ are called neighboring vertices or $(x,y)$ is an edge. The set of edges is denoted by $E$. The quadruple $(V,E, \omega ,\mu)$ is called a weighted graph (from now on, we only  write $V$ for short).

 We say that $V$ is locally finite if for any $x\in V$, there are finitely many neighboring vertices of $x$.  
A sequence $\left(x_k\right)_{k=1}^n$ of vertices  is  a path if $x_k \sim x_{k+1} $ for all $ k=1, \ldots, n-1$. A graph $V$ is called connected if, for any two vertices $x, y \in V$, there exists a path joining $x$ to $y$.

Throughout this paper, we always assume that $V$ is a locally finite, connected and weighted graph.
Let $\ell(V)$ be the set of  real-valued functions on $V$. The Laplace operator $\Delta$ is  defined on  $\ell(V)$ by 
$$\Delta f(x):=\sum_{y\sim x}\frac{\omega_{xy}}{\mu(x)}(f(y)-f(x))\mbox{ for }x\in V.$$
In this paper, we are interested in a system of  wave inequalities
\begin{equation}\label{e48231b}
\begin{cases}u_{t t}-\Delta u \geq  h_1|u|^p & \text { in } V \times(0, \infty), \\ 
	v_{t t}-\Delta v \geq h_2|v|^q& \text { in } V \times(0, \infty), \\ 
	u=u_0;\;v=v_0 & \text { in } V \times\{0\}, \\ u_t=u_1;\;v_t=v_1  & \text { in } V \times\{0\},
	\end{cases}
\end{equation}
where $p,q>1$, $h_1,h_2$ are positive potentials.

Next, we would like to present some assumptions on graphs, see  \cite{MPS24,MPS23}. Let us recall that 
a pseudo-metric on $V$ is a symmetric map $d: V \times V \rightarrow[0,\infty)$ satisfying $d(x,x)=0$ and the triangle inequality
$$
d(x, y) \leqslant d(x, z)+d(z, y)\mbox{ for  all }x, y, z \in V .$$
Given a  pseudo-metric $d$, define the jump size $j>0$ by
$$
j:=\sup \{d(x, y)~;~ x, y \in V, x\sim y\} .
$$
For any $x\in V$, $R>0$, the ball centered at $x$ and of radius $R>0$ related to the pseudo-metric $d$ is defined by
$$B_R(x)=\{y\in V~;~d(x,y)\le R\}.$$

\noindent{\textbf {Assumption A}} 
{\it \begin{enumerate}
		\item[(1)] there exists a constant $C_1>0$ such that $\sum\limits_{y\sim x}\omega_{xy}\le C_1\mu(x)$ for every $x\in V$;
		\item[(2)] there exists a pseudo-metric $d$ such that the jump size $j$ is finite;
		\item[(3)] the ball $B_R(x)$ with respect to $d$ is a finite set for any $x\in V$ and $ R>0$;
		\item[(4)] there exist  $x_0\in V$, $R_0>0,\alpha\in[0,1]$ and a constant $C_2>0$ such that
		$$\Delta d(x,x_0)\le\dfrac{C_2}{d^\alpha(x,x_0)}
	\mbox{ for any  }x\in V\setminus B_{R_0}(x_0).$$
\end{enumerate}}
Under  Assumption A, the authors in \cite{MPS24} obtained the following.

\noindent{\bf Theorem A. }{\it  Let  Assumption A be satisfied. Assume that $p>1$ and $\theta_1 \geqslant 2, \theta_2 \geqslant 2$ such that $\frac{2\theta_1}{\theta_2} \geqslant 1+\alpha$. Put 
	\begin{equation}\label{eer}
		E_R:=\left\{(x, t) \in V \times[0,\infty); R^{\theta_1} \leqslant d\left( x,x_0\right)^{\theta_1}+t^{\theta_2} \leqslant 2 R^{\theta_1}\right\}.
	\end{equation}
	Suppose that for every $R \geqslant R_0>0$,
\begin{equation}\label{eds106241}
		\int_0^{\infty} \sum_{x \in V} h^{-\frac{1}{p-1}} \mathds{1}_{E_R}(x, t) \mu(x) d t \leqslant C R^{\frac{(1+\alpha) p}{p-1}}
\end{equation}
and 
$$\liminf _{R \rightarrow \infty}\left\{\sum_{x \in B_R\left(x_0\right)} u_1^{+}(x) \mu(x)-\sum_{x \in B_{2 R}\left(x_0\right)} u_1^{-}(x) \mu(x)\right\} \geqslant 0.$$
 Then any non-negative weak solution of the inequality \eqref{eat3} in $V\times [0,\infty)$ must be trivial.}

As defined in \cite{MPS26}, a  weak solution $u$ of \eqref{eat3} means 
 $$u(x, \cdot) \in L_{l o c}^1([0, \infty)) \cap L_{l o c}^p([0, \infty), h(x, t) d t)\mbox{ for every } x \in V,$$  and
$$
\begin{aligned}
&	\int_0^{\infty} \sum_{x \in V} u(x, t) \varphi_{t t}(x, t) \mu(x) d t  -\int_0^{\infty} \sum_{x \in V} \Delta u(x, t) \varphi(x, t) \mu(x) d t\\
&+\sum_{x \in V} u_0(x) \varphi_t(x, 0) \mu(x) -\sum_{x \in V} u_1(x) \varphi(x, 0) \mu(x) \\
&\geqslant \int_0^{\infty} \sum_{x \in V} h(x, t)|u(x, t)|^p \varphi(x, t) \mu(x) d t
\end{aligned}
$$
for every compactly supported function $\varphi: V \times[0, \infty) \rightarrow \mathbb{R}$ such that $\varphi \geqslant 0$,  $\varphi(x, \cdot) \in C^2([0, \infty))$.

Similarly, we define a weak solution $(u,v)$ of \eqref{e48231b} as:
 $$u(x, \cdot) \in L_{l o c}^1([0, \infty)) \cap L_{l o c}^q([0, \infty), h_2(x, t) d t)\mbox{ for every } x \in V,$$  
$$v(x, \cdot) \in L_{l o c}^1([0, \infty)) \cap L_{l o c}^p([0, \infty), h_1(x, t) d t)\mbox{ for every } x \in V,$$ 
and 
$$
\begin{aligned}
	&	\int_0^{\infty} \sum_{x \in V} u(x, t) \varphi_{t t}(x, t) \mu(x) d t  -\int_0^{\infty} \sum_{x \in V} \Delta u(x, t) \varphi(x, t) \mu(x) d t\\
	&+\sum_{x \in V} u_0(x) \varphi_t(x, 0) \mu(x) -\sum_{x \in V} u_1(x) \varphi(x, 0) \mu(x) \\
	&\geqslant \int_0^{\infty} \sum_{x \in V} h_1(x, t)|v(x, t)|^p \varphi(x, t) \mu(x) d t,
\end{aligned}
$$
$$
\begin{aligned}
	&	\int_0^{\infty} \sum_{x \in V} v(x, t) \varphi_{t t}(x, t) \mu(x) d t  -\int_0^{\infty} \sum_{x \in V} \Delta v(x, t) \varphi(x, t) \mu(x) d t\\
	&+\sum_{x \in V} v_0(x) \varphi_t(x, 0) \mu(x) -\sum_{x \in V} v_1(x) \varphi(x, 0) \mu(x) \\
	&\geqslant \int_0^{\infty} \sum_{x \in V} h_2(x, t)|u(x, t)|^q \varphi(x, t) \mu(x) d t
\end{aligned}
$$
for every compactly supported function $\varphi: V \times[0, \infty) \rightarrow \mathbb{R}$ such that $\varphi \geqslant 0$,  $\varphi(x, \cdot) \in C^2([0, \infty))$.

Our first main result in this paper is the following.
\begin{theorem}\label{th5}
	Let  Assumption A be satisfied. Assume that $p\geq  q>1$,  and $\theta_1 \geqslant 2, \theta_2 \geqslant 2$ such that $\frac{2\theta_1}{\theta_2} \geqslant 1+\alpha$. Suppose that for every $R \geqslant R_0>0$,
	\begin{align}\label{eds411b}
\begin{split}
		&	\int_{0}^\infty \sum_{x \in V} \mathds{1}_{E_R}(x, t) h_1^{-\frac{1}{p-1}}\mu(x) d t \leqslant C R^{\frac{p(q+1)(1+\alpha)}{pq-1}},\\
	&	\int_{0}^\infty \sum_{x \in V} \mathds{1}_{E_R}(x, t) h_2^{-\frac{1}{q-1}}\mu(x) d t \leqslant C R^{\frac{p(q+1)(1+\alpha)}{pq-1}},
\end{split}
	\end{align}
	where
	$E_R$ is defined in \eqref{eer}. Then, any non-negative weak solution  of the system of  \eqref{e48231b}  must be trivial.
\end{theorem}
Remark that, when $p=q$ and $h_1=h_2=h$, the condition \eqref{eds411b} becomes \eqref{eds106241}. Therefore, our result is a natural  extension of the result in \cite{MPS26} from a wave  inequality to a system of wave  inequalities.  

A direct consequence of Theorem \ref{th5} when $h_1=h_2=1$ is as follows.
\begin{corollary}
		Let  Assumption A be satisfied. Assume that $p\geq  q>1$,  and $\theta_1 \geqslant 2, \theta_2 \geqslant 2$ such that $\frac{2\theta_1}{\theta_2} \geqslant 1+\alpha$. Suppose that for every $R \geqslant R_0>0$,
	\begin{align*}
		\begin{split}
			&	\int_{0}^\infty \sum_{x \in V} \mathds{1}_{E_R}(x, t) \mu(x) d t \leqslant C R^{\frac{p(q+1)(1+\alpha)}{pq-1}},\\
		\end{split}
	\end{align*}
	where
	$E_R$ is defined in \eqref{eer}. Then, any non-negative weak solution  of the system 
	$$\begin{cases}u_{t t}-\Delta u \geq  |u|^p & \text { in } V \times(0, \infty), \\ 
		v_{t t}-\Delta v \geq |v|^q& \text { in } V \times(0, \infty), \\ 
		u=u_0;\;v=v_0 & \text { in } V \times\{0\}, \\ u_t=u_1;\;v_t=v_1  & \text { in } V \times\{0\},
	\end{cases}$$
	 must be trivial.
\end{corollary}
Without any assumption on the sign of solutions, as in \cite{MPS24}, under the assumption A, one defines the space
$$
X_\delta=\left\{f\in \ell(V);\sum_{x \in V}f(x) e^{-\delta d\left(x, x_0\right)} \mu(x)<\infty\right\},
$$
where $\delta>0$.

Our second result in this paper reads as follows.
\begin{theorem}\label{th2}
Assume that Assumption A holds and  $p\geq q>1$. Let $h_1,h_2$ be  positive functions. Suppose that $u$ is a weak solution of \eqref{e48231b} with 
$$
\sum_{x \in V} u_1(x) \mu(x) \geqslant 0, \; \sum_{x \in V} v_1(x) \mu(x) \geqslant 0,
$$
$u_0, u_1,v_0,v_1 \in X_\delta$, and $u,v \in L_{\text {loc }}^1\left([0, \infty), X_\delta\right)$, $\delta>0$. Assume, in addition that, for every $R \geqslant R_0>1$, 
$$\max(V_1(R), V_2(R),V_3(R), V_4(R)) \leqslant C R^{\frac{p(q+1)(1+\alpha)}{pq-1}},$$
where
$$
V_1(R)=\int_{R^{\frac{1+\alpha}{2}}}^{2 R^{\frac{1+\alpha}{2}}} \sum_{x \in B_R\left(x_0\right)} h_1^{-\frac{1}{p-1}}(x, t) e^{-\delta \frac{d\left(x, x_0\right)}{R}} \mu(x) d t, 
$$
$$
V_2(R)=\int_0^{2 R^{\frac{1+\alpha}{2}}} \sum_{x \in V \backslash B_R\left(x_0\right)} h_1^{-\frac{1}{p-1}}(x, t) e^{-\delta \frac{d\left(x, x_0\right)}{R}} \mu(x) d t 
$$
$$V_3(R)=\int_{R^{\frac{1+\alpha}{2}}}^{2 R^{\frac{1+\alpha}{2}}} \sum_{x \in B_R\left(x_0\right)} h_2^{-\frac{1}{q-1}}(x, t) e^{-\delta \frac{d\left(x, x_0\right)}{R}} \mu(x) d t,
$$
and 
$$
V_4(R)=\int_0^{2 R^{\frac{1+\alpha}{2}}} \sum_{x \in V \backslash B_R\left(x_0\right)} h_2^{-\frac{1}{q-1}}(x, t) e^{-\delta \frac{d\left(x, x_0\right)}{R}} \mu(x) d t,
$$
then $u=v= 0$.
\end{theorem}
This result can be seen as a generalization of \cite[Theorem 3.4]{MPS26} from a single inequality to a system of inequalities. On the other hand, in the case $\alpha=1$, our exponent in the above theorems  coincides with that found in Euclidean setting, \cite{JMY20}. In particular, when $p=q$, it goes back to that in \cite{Ka80}. Recall that these exponents  were  shown to be sharp  in Euclidean setting.

\section{Proof of Theorem \ref{th5}}\label{s4}

 Let us first recall an important result established in \cite{MPS26}. Denote by  $\varphi \in C^2([0,\infty))$  a cut-off function such that $\varphi = 1$ in $[0,1]$, $\varphi = 0$ in $[2,\infty)$, and $\varphi^{\prime} \leqslant 0$. For each $x \in V, t \in[0,\infty)$ and $R \geqslant R_0$, put
$$
 \varphi_R(x, t)=\varphi\left(\frac{t^{\theta_2}+d\left(x_0, x\right)^{\theta_1}}{R^{\theta_1}}\right).
$$
Then we have the following lemma which was given in \cite{MPS26}.

\begin{lemma}\label{l-MPS}Let 
	$$F_R:=\left\{(x, t) \in V \times[0,\infty);(R / 2)^{\theta_1} \leqslant d\left(x_0, x\right)^{\theta_1}+t^{\theta_2} \leqslant(4 R)^{\theta_1}\right\}.$$
	Then, there exists a positive constant $C$ such that
$$
-\Delta \varphi_R(x, t) \leqslant \frac{C}{R^{1+\alpha}} \mathds{1}_{F_R}(x, t), \mbox{ for all }(x,t)\in V\times[0,\infty), 
$$
and
$$
	- (\varphi_R)_t(x, t) \leqslant \frac{C}{R^{\frac{\theta_1}{\theta_2}}} \mathds{1}_{E_R}(x, t)\mbox{ for all }(x,t)\in V\times[0,\infty).
$$
and
$$
- (\varphi_R)_{tt}(x, t) \leqslant \frac{C}{R^{2\frac{\theta_1}{\theta_2}}} \mathds{1}_{E_R}(x, t)\mbox{ for all }(x,t)\in V\times[0,\infty).
$$
\end{lemma}
Let $s$ be a positive parameter which will be chosen sufficiently large. It follows from \eqref{e48231b} that 
\begin{align}\label{eds266241}
\begin{split}
	&	\int_0^{\infty} \sum_{x \in V} \mu(x) h_1(x,t)v^p(x, t) \varphi_R^s(x, t) d t\\
		& \leqslant  \int_0^{\infty} \sum_{x \in V} \mu(x) u(x, t) (\varphi_R^{s})_{tt}(x, t)  d t \\
	& -\int_0^{\infty} \sum_{x \in V} \mu(x) \Delta u(x, t) \varphi_R^s(x, t) d t  -\sum_{x \in V} \mu(x) u_t(x, 0) \varphi_R^s(x, 0) \\
	&+\sum_{x \in V} \mu(x) u(x, 0) (\varphi_R^s)_t(x, 0).
\end{split}
\end{align}
Note that $(\varphi_R^s)_t(x, 0)=s\varphi_R^{s-1}(\varphi_R)_t(x,0)=0$ since $\theta_2\geq 2$. Then the last term in \eqref{eds266241} vanishes. On the other hand, 
$$|(\varphi_R^{s})_{tt}|=|s(s-1)\varphi_R^{s-2}(\varphi_R)^2_t+s\varphi_R^{s-1}(\varphi_R)_{tt}|\leq \frac{C}{R^{2\frac{\theta_1}{\theta_2}}}\varphi_R^{s-2}\mathds{1}_{E_R}(x, t).$$
Therefore, \eqref{eds266241} implies
\begin{align}\label{eds266241b}
	\begin{split}
			&	\int_0^{\infty} \sum_{x \in V} \mu(x) h_1(x,t)v^p(x, t) \varphi_R^s(x, t) d t\\
	&\leqslant  \frac{C}{R^{2\frac{\theta_1}{\theta_2}}}\int_0^{\infty} \sum_{x \in V} \mu(x) u(x, t) \varphi_R^{s-2}(x, t) \mathds{1}_{E_R}(x, t)d t \\
	& -\int_0^{\infty} \sum_{x \in V} \mu(x) \Delta u(x, t) \varphi_R^s(x, t) d t-\sum_{x \in V} \mu(x) u_1(x) \varphi_R^s(x, 0).
\end{split}
\end{align}
Since $s>1$, we have
\begin{align}\label{eds266242}
	\begin{split}
		-\Delta \varphi_R^s(x, t) & =-\frac{1}{\mu(x)} \sum_{y \sim x} \omega_{x y}\left(\varphi_R^s(y, t)-\varphi_R^s(x, t)\right) \\
		& \leqslant-\frac{1}{\mu(x)} \sum_{y \sim x} s \omega_{x y} \varphi_R^{s-1}(x, t)\left(\varphi_R(y, t)-\varphi_R(x, t)\right) \\
		& =-s \varphi_R^{s-1}(x, t) \Delta \varphi_R(x, t).
	\end{split}
\end{align}
By the definition of $\varphi_R$, we get
\begin{align}\label{eds266242b}
	\begin{split}
	&	\sum_{x \in V} \mu(x) u_1(x) \varphi_R^s(x, 0)=\sum_{x \in B_{2R}(x_0)} \mu(x) (u_1^+(x)-u_1^-(x)) \varphi_R^s(x, 0)\\
	&\geq \sum_{x \in B_{R}(x_0)} \mu(x) u_1^+(x))-\sum_{x \in B_{2R}(x_0)} \mu(x) u_1^+(x),
	\end{split}
\end{align}
where $u^+=\frac{|u|+u}{2}$ and $u^-=\frac{|u|-u}{2}$.
 
Substituting  \eqref{eds266242}, \eqref{eds266242b} into \eqref{eds266241b} and  using Lemma \ref{l-MPS}, we obtain
\begin{align}\label{eds266245}
	\begin{split}
		&	\int_0^{\infty} \sum_{x \in V} \mu(x) h_1(x,t)v^p(x, t) \varphi_R^s(x, t) d t\\
		&\leqslant  \frac{C}{R^{2\frac{\theta_1}{\theta_2}}}\int_0^{\infty} \sum_{x \in V} \mu(x) u(x, t) \varphi_R^{s-2}(x, t) \mathds{1}_{E_R}(x, t)d t \\
		& +\frac{C}{R^{1+\alpha}} \int_0^{\infty} \sum_{x \in V} \mu(x) u(x, t) \varphi_R^{s-1}(x, t) \mathds{1}_{F_R}(x, t) d t\\
		&-\sum_{x \in B_{R}(x_0)} \mu(x) u_1^+(x)+\sum_{x \in B_{2R}(x_0)} \mu(x) u_1^+(x)\\
		&\leq \frac{C}{R^{\alpha+1}} \int_0^{\infty} \sum_{x \in V} \mu(x) u(x, t) \varphi_R^{s-2}(x, t) \mathds{1}_{F_R}(x, t) d t\\
		&-\sum_{x \in B_{R}(x_0)} \mu(x) u_1^+(x)+\sum_{x \in B_{2R}(x_0)} \mu(x) u_1^+(x),
	\end{split}
\end{align}
where we have used $E_R\subset F_R $, $0\leq \varphi_R\leq 1$ and $\frac{2\theta_1}{\theta_2}\geq \alpha+1$.
Using the H\"{o}lder inequality, we obtain 
\begin{align}\label{eds266246}
\begin{split}
	&	\int_0^{\infty} \sum_{x \in V} \mu(x) u(x, t) \varphi_R^{s-2}(x, t) \mathds{1}_{F_R}(x, t) d t\\
	&\leq \left(	\int_0^{\infty} \sum_{x \in V} \mu(x) h_2(x,t)u^q(x, t) \varphi_R^{q(s-2)}(x, t) \mathds{1}_{F_R}(x, t)d t\right)^\frac{1}{q}\times\\
	&\times\left(  \int_0^{\infty} \sum_{x \in V} \mu(x)h_2^{-\frac{1}{q-1}}(x,t)\mathds{1}_{F_R}(x, t)dt\right)^{\frac{q-1}{q}}.
\end{split}
\end{align}
By using similar argument as in the proof of  \eqref{eds266245} for the second inequality in \eqref{e48231b}, we also have 
\begin{equation}\label{eds266247}
	\begin{aligned}
		\int_0^{\infty} \sum_{x \in V} \mu(x) & h_2(x,t) u^q(x, t) \varphi_R^{q(s-2)}(x, t) d t \\
		& \leqslant \frac{C}{R^{\alpha+1}} \int_0^{\infty} \sum_{x \in V} \mu(x) v(x, t) \varphi_R^{q(s-2)-2}(x, t) \mathds{1}_{F_R}(x, t) d t\\
		&-\sum_{x \in B_{R}(x_0)} \mu(x) v_1^+(x)+\sum_{x \in B_{2R}(x_0)} \mu(x) v_1^+(x),
	\end{aligned}
\end{equation}
Applying again the H\"{o}lder inequality to the right hand side of \eqref{eds266247}, we arrive at
\begin{align}\label{eds266248}
	\begin{split}
		&	\int_0^{\infty} \sum_{x \in V} \mu(x) v(x, t) \varphi_R^{q(s-2)-2}(x, t) \mathds{1}_{F_R}(x, t) d t\\
		&\leq \left(	\int_0^{\infty} \sum_{x \in V} \mu(x) h_1(x,t)v^p(x, t) \varphi_R^{p(q(s-2)-2)}(x, t) \mathds{1}_{F_R}(x, t)d t\right)^\frac{1}{p}\times\\
		&\times\left(  \int_0^{\infty} \sum_{x \in V} \mu(x)h_1^{-\frac{1}{p-1}}(x,t)\mathds{1}_{F_R}(x, t)dt\right)^{\frac{p-1}{p}}\\
			&\leq \left(	\int_0^{\infty} \sum_{x \in V} \mu(x) h_1(x,t)v^p(x, t) \varphi_R^{s}(x, t) \mathds{1}_{F_R}(x, t)d t\right)^\frac{1}{p}\times \\
			&\times\left(  \int_0^{\infty} \sum_{x \in V} \mu(x)h_1^{-\frac{1}{p-1}}(x,t)\mathds{1}_{F_R}(x, t)dt\right)^{\frac{p-1}{p}},
	\end{split}
\end{align}
where  $s$ is chosen sufficiently large such that $p(q(s-2)-2)>s$.

The assumption on the initial data $u_1,v_1$ implies that there exists a sequence $(R_n), R_n\to \infty$ such that 
$$-\sum_{x \in B_{R_n}(x_0)} \mu(x) u_1^+(x)+\sum_{x \in B_{2R_n}(x_0)} \mu(x) u_1^+(x)\leq 0$$
and 
$$-\sum_{x \in B_{R_n}(x_0)} \mu(x) v_1^+(x)+\sum_{x \in B_{2R_n}(x_0)} \mu(x) v_1^+(x)\leq 0.$$
From this and  \eqref{eds266245}, \eqref{eds266246}, \eqref{eds266247}, \eqref{eds266248}, we arrive at 
\begin{equation}\label{eds266249b}
	\begin{aligned}
		&	\int_0^{\infty} \sum_{x \in V} \mu(x) h_1(x,t) v^p(x, t) \varphi_{R_n}^s(x, t) d t \\
		&\leq  \frac{C}{{R_n}^{(\alpha+1)(1+\frac{1}{q})}} \left(\int_0^{\infty} \sum_{x \in V} \mu(x)  h_1(x,t)v^p(x, t) \varphi_{R_n}^s(x, t) \mathds{1}_{F_{R_n}}(x, t)d t\right)^{\frac{1}{pq}}\times\\
		&\times\left(  \int_0^{\infty} \sum_{x \in V} \mu(x)h_2^{-\frac{1}{q-1}}(x,t)\mathds{1}_{F_{R_n}}(x, t)dt\right)^{\frac{q-1}{q}}\times\\
	&	\times\left(  \int_0^{\infty} \sum_{x \in V} \mu(x)h_1^{-\frac{1}{p-1}}(x,t)\mathds{1}_{F_{R_n}}(x, t)dt\right)^{\frac{p-1}{pq}}.
	\end{aligned}
\end{equation}
This follows that 
\begin{equation}\label{eds266249}
	\begin{aligned}
&	\left(\int_0^{\infty} \sum_{x \in V} \mu(x)h_1(x,t)  v^p(x, t) \varphi_R^s(x, t) d t \right)^{1-\frac{1}{pq}}\\
	 &\leq  \frac{C}{R^{(\alpha+1)(1+\frac{1}{q})}}\left(  \int_0^{\infty} \sum_{x \in V} \mu(x)h_2^{-\frac{1}{q-1}}(x,t)\mathds{1}_{F_{R_n}}(x, t)dt\right)^{\frac{q-1}{q}}\\
	 &\times\left(  \int_0^{\infty} \sum_{x \in V} \mu(x)h_1^{-\frac{1}{p-1}}(x,t)\mathds{1}_{F_{R_n}}(x, t)dt\right)^{\frac{p-1}{pq}}.	\end{aligned}
\end{equation}
Let $m$ be a fixed positive integer large enough such that $m>3\theta_1-1$,{ and set
	$\kappa=\kappa(k):=2^{\frac{k}{\theta_1}-1}$. Then one has
	$F_{R_n}\subset \cup_{k=0}^mE_{\kappa {R_n}}.$
	This implies  
	$$
	\int_0^{\infty} \sum_{x \in V} \mu(x)h_1^{-\frac{1}{p-1}}(x,t)\mathds{1}_{F_{R_n}}(x, t)dt  \leq \sum_{k=0}^m\int_0^{\infty} 
	\sum_{x \in V} \mu(x)h_1^{-\frac{1}{p-1}}(x,t)\mathds{1}_{E_{\kappa {R_n}}}(x, t)dt
	$$
	and 
		$$
	\int_0^{\infty} \sum_{x \in V} \mu(x)h_2^{-\frac{1}{q-1}}(x,t)\mathds{1}_{F_{R_n}}(x, t)dt  \leq \sum_{k=0}^m\int_0^{\infty} 
	\sum_{x \in V} \mu(x)h_2^{-\frac{1}{q-1}}(x,t)\mathds{1}_{E_{\kappa {R_n}}}(x, t)dt
	$$
	Now we use the assumption \eqref{eds411b}, one has
	$$
	\int_0^{\infty} 
	\sum_{x \in V} \mu(x)h_2^{-\frac{1}{q-1}}(x,t)\mathds{1}_{F_{R_n}}(x, t)dt\le C
	{R_n}^{\frac{p(q+1)(1+\alpha)}{pq-1}}
	$$
	and 
		$$
	\int_0^{\infty} 
	\sum_{x \in V} \mu(x)h_1^{-\frac{1}{p-1}}(x,t)\mathds{1}_{F_{R_n}}(x, t)dt\le C
		{R_n}^{\frac{p(q+1)(1+\alpha)}{pq-1}}.
		$$
Combining this with \eqref{eds266249}, we deduce that 
	\begin{equation}\label{eds2662410}
		\left(	\int_0^{\infty} \sum_{x \in V} \mu(x)h_1(x,t)  v^p(x, t) \varphi_{R_n}^s(x, t) d t\right)^\frac{pq-1}{pq}  \leq C.
	\end{equation}
}
Letting $n\to\infty$, i.e $R_n\to\infty$ in \eqref{eds2662410}, we get
$$\int_0^{\infty} \sum_{x \in V} \mu(x) h_1(x,t) v^p(x, t)  d t<\infty.$$
This inequality and \eqref{eds266249b} imply that 
\begin{equation}\label{eds2662411}
	\begin{aligned}
		&	\int_0^{\infty} \sum_{x \in V} \mu(x) h_1(x,t) v^p(x, t) \varphi_{R_n}^s(x, t) d t \\
		&\leq  C \left(\int_0^{\infty} \sum_{x \in V} \mu(x)h_1(x,t)  v^p(x, t) \varphi_{R_n}^s(x, t) \mathds{1}_{F_{R_n}}(x, t)d t\right)^{\frac{1}{pq}}\to 0 \mbox{ as }n\to\infty.
	\end{aligned}
\end{equation}
Since $v$ is non-negative, we obtain $v=0$. Similarly, we also arrive at $u=0$. The proof is completed.\qed
\section{Proof of Theorem \ref{th2}}
We first recall the construction of test functions introduced in \cite{MPS24}.

\noindent{\bf Beginning of the proof}

Let $j$  be the jump side given in Assumption A. Let $\psi $  be a $C^2$ positive function defined on $[-j, \infty)$ satisfying $\psi^{\prime} \leqslant 0$,  $\psi \equiv 1$ on $[-j, 1]$ and  $\psi(r)=e^{-\delta r}$ when  $r \geqslant 2$. Here $\delta$ is given in the definition of the space $X_\delta$.

Let $\eta:[0, \infty) \rightarrow[0, \infty)$ be a $C^2$ compactly supported function satisfying  $\eta^{\prime} \leqslant 0$, 
$\eta \equiv 1$ on $[0,1]$ and $\eta \equiv 0$ on $[2, \infty)$.

Let $R>0$ be a sufficiently large parameter, define  next the test function on $V\times [0,\infty)$:
$$
\varphi_R(x, t)=\eta^s\left(\frac{t}{R^{\frac{1+\alpha}{2}}}\right) \psi\left(\frac{d\left(x, x_0\right)-j}{R}\right),
$$
where $s$ is large enough chosen later on.  Then, the following lemma was proved in \cite{MPS24}.
\begin{lemma}\label{l3} Let  $Q_R:=V \times\left[R^{\frac{1+\alpha}{2}}, 2 R^{\frac{1+\alpha}{2}}\right]$. Then, we have  $ $
		$$
		\left|\left(\varphi_R\right)_t(x, t)\right| \leqslant \frac{C}{R^{\frac{1+\alpha}{2}}} \eta^{s-1}\left(\frac{t}{R^{\frac{1+\alpha}{2}}}\right) e^{-\delta \frac{d\left(x, x_0\right)}{R}} \mathds{1}_{Q_R}(x, t),
		$$
			$$
			\left|\left(\varphi_R\right)_{t t}(x, t)\right| \leqslant \frac{C}{R^{1+\alpha}} \eta^{s-2}\left(\frac{t}{R^{\frac{1+\alpha}{2}}}\right) e^{-\delta \frac{d\left(x, x_0\right)}{R}} \mathds{1}_{Q_R}(x, t) ,
			$$
and 
			$$
			\left|\Delta \varphi_R(x, t)\right| \leqslant \frac{C}{R^{1+\alpha}} \eta^s\left(\frac{t}{R^{\frac{1+\alpha}{2}}}\right) e^{-\delta \frac{d\left(x, x_0\right)}{R}} \mathds{1}_{V \backslash B_R\left(x_0\right)}(x) .
			$$

\end{lemma}
We now prove Theorem \ref{th2}. By definition of weak solution of \eqref{e48231b}, we have 
\begin{equation}\label{eat41}
	\begin{aligned} &\int_0^{\infty} \sum_{x \in V} \mu(x) h_1(x,t)|v(x, t)|^p \varphi_R(x, t) d t \leqslant  -\int_0^{\infty} \sum_{x \in V} \Delta u(x, t) \varphi_R(x, t) \mu(x) d t \\ & +\sum_{x \in V} u_0(x)\left(\varphi_R\right)_t(x, 0) \mu(x)-\sum_{x \in V} u_1(x) \varphi_R(x, 0) \mu(x) \\ & +\int_0^{\infty} \sum_{x \in V} u(x, t)\left(\varphi_R\right)_{t t}(x, t) \mu(x) d t.\end{aligned}
\end{equation}
For the last term in the right hand side of this inequality, using Lemma \ref{l3}, we have 
\begin{align}\label{ew2}
	 \begin{split}
		& \left|\int_0^{\infty} \sum_{x \in V} u(x, t)\left(\varphi_R\right)_{t t}(x, t) \mu(x) d t\right|\leqslant \int_0^{\infty} \sum_{x \in V}| u(x, t) |\left(\varphi_R\right)_{t t}(x, t) \mid \mu(x) d t \\ \leqslant & \frac{C}{R^{1+\alpha}} \int_{R^{\frac{1+\alpha}{2}}}^{2 R^{\frac{1+\alpha}{2}}} \sum_{x \in V}|u(x, t)| \eta^{s-2}\left(\frac{t}{R^{\frac{1+\alpha}{2}}}\right) e^{-\delta \frac{d\left(x, x_0\right)}{R}} \mu(x) d t \\
	 \leqslant &  \frac{C}{R^{1+\alpha}} \left(\int_{R^{\frac{1+\alpha}{2}}}^{2 R^{\frac{1+\alpha}{2}}} \sum_{x \in V}|u(x, t)|^q h_2(x, t) \eta^s\left(\frac{t}{R^{\frac{1+\alpha}{2}}}\right) e^{-\delta \frac{d\left(x, x_0\right)}{R}} \mu(x) d t \right)^\frac{1}{q}\\ 
	 &\times \left( \int_{R^{\frac{1+\alpha}{2}}}^{2 R^{\frac{1+\alpha}{2}}} \sum_{x \in V} h_2^{-\frac{1}{q-1}}(x, t) \eta^{s-\frac{2 q}{q-1}}\left(\frac{t}{R^{\frac{1+\alpha}{2}}}\right) e^{-\delta \frac{d\left(x, x_0\right)}{R}} \mu(x) d t\right)^\frac{q-1}{q}.
\end{split}\end{align}
Next, we estimate the first term of \eqref{eat41} by using again Lemma \ref{l3}
\begin{align}\label{ew3}
	\begin{split}
		 &\left|- \int_0^{\infty} \sum_{x \in V} \Delta u(x, t) \varphi_R(x, t) \mu(x) d t \right| \\ &=  \left|\int_0^{\infty} \sum_{x \in V} u(x, t) \Delta \varphi_R(x, t) \mu(x) d t\right| \\ & \leqslant \int_0^{\infty} \sum_{x \in V}|u(x, t)|\left|\Delta \varphi_R\right| \mu(x) d t \\ & \leqslant \frac{C}{R^{1+\alpha}} \int_0^{\infty} \sum_{x \in V \backslash B_R\left(x_0\right)}|u(x, t)| \eta^s\left(\frac{t}{R^{\frac{1+\alpha}{2}}}\right) e^{-\delta \frac{d\left(x, x_0\right)}{R}} \mu(x) d t \\ & \leqslant \frac{C}{R^{1+\alpha}} \left( \int_0^{\infty} \sum_{x \in V \backslash B_R\left(x_0\right)}|u(x, t)|^q h_2(x, t) \eta^s\left(\frac{t}{R^{\frac{1+\alpha}{2}}}\right) e^{-\delta \frac{d\left(x, x_0\right)}{R}} \mu(x) d t\right)^\frac{1}{q} \\ & \times \left(\int_0^{\infty} \sum_{x \in V \backslash B_R\left(x_0\right)} h_2^{-\frac{1}{q-1}}(x, t) \eta^s\left(\frac{t}{R^{\frac{1+\alpha}{2}}}\right) e^{-\delta \frac{d\left(x, x_0\right)}{R}} \mu(x) d t\right)^{\frac{q-1}{q}}.
	\end{split}
\end{align}
For the second term in \eqref{eat41}, we have
$$
\left(\varphi_R\right)_t(x, t)=\frac{s}{R^{\frac{1+\alpha}{2}}} \eta^{s-1}\left(\frac{t}{R^{\frac{1+\alpha}{2}}}\right) \psi\left(\frac{d\left(x, x_0\right)-j}{R}\right) \eta^{\prime}\left(\frac{t}{R^{\frac{1+\alpha}{2}}}\right)\mathds{1}_{Q_R}(x, t).
$$
Hence, $\left(\varphi_R\right)_t(x, 0)=0$, which implies

\begin{equation}\label{ew4}
	\sum_{x \in V} u_0(x)\left(\varphi_R\right)_t(x, 0) \mu(x)=0 .
\end{equation}

Finally, there holds

\begin{equation}\label{ew8}
	\begin{aligned}
	-\sum_{x \in V} u_1(x) \varphi_R(x, 0) \mu(x)= & -\sum_{x \in V} u_1(x) \mu(x) \psi\left(\frac{d\left(x, x_0\right)-j}{R}\right) \\
	= & -\sum_{x \in V} u_1^{+}(x) \mu(x) \psi\left(\frac{d\left(x, x_0\right)-j}{R}\right) \\
	& +\sum_{x \in V} u_1^{-}(x) \mu(x) \psi\left(\frac{d\left(x, x_0\right)-j}{R}\right) \\
	\leqslant & -\sum_{x \in B_{R+j}\left(x_0\right)} u_1^{+}(x) \mu(x)+\sum_{x \in V} u_1^{-}(x) \mu(x) .
\end{aligned}
\end{equation}
Combining this and \eqref{ew2}, \eqref{ew3}, \eqref{ew4} we have
\begin{align}\label{eat11}
\begin{split}
	&	\int_0^{\infty} \sum_{x \in V} \mu(x) h_1(x,t)|v(x, t)|^p \varphi_R(x, t) d t\\
&\leq \frac{C}{R^{1+\alpha}} \left( \int_0^{\infty} \sum_{x \in V }|u(x, t)|^q h_2(x, t) \eta^s\left(\frac{t}{R^{\frac{1+\alpha}{2}}}\right) e^{-\delta \frac{d\left(x, x_0\right)}{R}} \mu(x) d t\right)^\frac{1}{q}\left(S_1+S_2\right)\\
	& -\sum_{x \in B_{R+j}\left(x_0\right)} u_1^{+}(x) \mu(x)+\sum_{x \in V} u_1^{-}(x) \mu(x),
\end{split}
\end{align}
where 
$$S_1=\left( \int_{R^{\frac{1+\alpha}{2}}}^{2 R^{\frac{1+\alpha}{2}}} \sum_{x \in V} h_2^{-\frac{1}{q-1}}(x, t) \eta^{s-\frac{2 q}{q-1}}\left(\frac{t}{R^{\frac{1+\alpha}{2}}}\right) e^{-\delta \frac{d\left(x, x_0\right)}{R}} \mu(x) d t\right)^\frac{q-1}{q}$$
and
$$S_2=\left(\int_0^{\infty} \sum_{x \in V \backslash B_R\left(x_0\right)} h_2^{-\frac{1}{q-1}}(x, t) \eta^s\left(\frac{t}{R^{\frac{1+\alpha}{2}}}\right) e^{-\delta \frac{d\left(x, x_0\right)}{R}} \mu(x) d t\right)^{\frac{q-1}{q}}.$$
Remark that, since $p\geq q$, we have 
\begin{align}\label{ew6}
\begin{split}
		R^{-(\alpha+1)}(S_1+S_2)&\leq CR^{-(\alpha+1)} R^{\frac{p(q+1)(\alpha+1)}{pq-1}\frac{q-1}{q}}\\
	&=CR^{(\alpha+1)\frac{-p+q}{(pq-1)q}}.
\end{split}
\end{align}
Similarly, by using the second inequality in \eqref{e48231b} we also obtain 
\begin{align}\label{eat21}
\begin{split}
	&	\int_0^{\infty} \sum_{x \in V \backslash B_R\left(x_0\right)}|u(x, t)|^q h_2(x, t) \varphi_R(x, t) \mu(x) d t\\
&\leq \frac{C}{R^{1+\alpha}} \left(	\int_0^{\infty} \sum_{x \in V} \mu(x) h_1(x,t)|v(x, t)|^p \varphi_R(x, t) d t\right)^\frac{1}{p}(S_3+S_4)\\
&-\sum_{x \in B_{R+j}\left(x_0\right)} v_1^{+}(x) \mu(x)+\sum_{x \in V} v_1^{-}(x) \mu(x),
\end{split}
\end{align}
where
$$S_3=\left( \int_{R^{\frac{1+\alpha}{2}}}^{2 R^{\frac{1+\alpha}{2}}} \sum_{x \in V} h_1^{-\frac{1}{p-1}}(x, t) \eta^{s-\frac{2 p}{p-1}}\left(\frac{t}{R^{\frac{1+\alpha}{2}}}\right) e^{-\delta \frac{d\left(x, x_0\right)}{R}} \mu(x) d t\right)^\frac{p-1}{p},$$
$$S_4=\left(\int_0^{\infty} \sum_{x \in V \backslash B_R\left(x_0\right)} h_1^{-\frac{1}{p-1}}(x, t) \eta^s\left(\frac{t}{R^{\frac{1+\alpha}{2}}}\right) e^{-\delta \frac{d\left(x, x_0\right)}{R}} \mu(x) d t\right)^{\frac{p-1}{p}}.$$
Notice again that, 
	\begin{align}\label{ew5}
		\begin{split}
	R^{-(\alpha+1)}(S_3+S_4)&\leq CR^{-(\alpha+1)} R^{\frac{p(q+1)(\alpha+1)}{pq-1}\frac{p-1}{p}}\\
	&=CR^{(\alpha+1)\frac{p-q}{pq-1}}.
\end{split}
\end{align}
On the other hand, by using the assumptions on the initial data, we have $$-\sum\limits_{x \in B_{R+j}\left(x_0\right)} u_1^{+}(x) \mu(x)+\sum\limits_{x \in V} u_1^{-}(x) \mu(x)$$  and $$-\sum\limits_{x \in B_{R+j}\left(x_0\right)} v_1^{+}(x) \mu(x)+\sum\limits_{x \in V} v_1^{-}(x) \mu(x)$$ are  decreasing and tend to zero as $R\to\infty$.
Thus, combining \eqref{eat11}, \eqref{ew6},  \eqref{eat21} and  \eqref{ew5} we arrive at 
\begin{equation*}
	\int_0^{\infty} \sum_{x \in V} \mu(x) h_1(x,t)|v(x, t)|^p d t\leq C
\end{equation*}
and
\begin{equation}\label{ew7}
	\int_0^{\infty} \sum_{x \in V} \mu(x) h_2(x,t)|u(x, t)|^q d t\leq C.
\end{equation}

\noindent{\bf End of the proof}

From this integrability, we will show the contradiction. On one hand, we have
\begin{equation}\label{eat51}
	\begin{aligned} & |\int_0^{\infty} \sum_{x \in V} \Delta u(x, t) \varphi_R(x, t) \mu(x) d t \mid \\ & \leqslant \int_0^{\infty} \sum_{x \in V}|u(x, t)|\left|\Delta \varphi_R\right| \mu(x) d t \\
	 \leqslant & \frac{C}{R^{1+\alpha}} \int_0^{\infty} \sum_{x \in V \backslash B_R\left(x_0\right)}|u(x, t)| \eta^s\left(\frac{t}{R^{\frac{1+\alpha}{2}}}\right) e^{-\delta \frac{d\left(x, x_0\right)}{R}} \mu(x) d t \\
	  \leqslant &\frac{C}{R^{1+\alpha}} \left(\int_0^{\infty} \sum_{x \in V \backslash B_R\left(x_0\right)}|u(x, t)|^q h_2(x, t) \eta^s\left(\frac{t}{R^{\frac{1+\alpha}{2}}}\right) e^{-\delta \frac{d\left(x, x_0\right)}{R}} \mu(x) d t\right)^{\frac{1}{q}} \\
	   & \times\left( \int_0^{\infty} \sum_{x \in V \backslash B_R\left(x_0\right)} h_2^{-\frac{1}{q-1}}(x, t) \eta^s\left(\frac{t}{R^{\frac{1+\alpha}{2}}}\right) e^{-\delta \frac{d\left(x, x_0\right)}{R}} \mu(x) d t\right)^{\frac{q-1}{q}} \\ 
	   \leqslant & C\left(\int_0^{\infty} \sum_{x \in V \backslash B_R\left(x_0\right)}|u(x, t)|^q h_2(x, t) \mu(x) d t\right)^{\frac{1}{q}}R^{(\alpha+1)\frac{-p+q}{(pq-1)q}}\to 0 \mbox{ as }R\to \infty,\end{aligned}
\end{equation}
where we have used \eqref{ew7}. 
	   In addition,
	   \begin{equation}\label{eat61}
	   	\begin{aligned} \mid \int_0^{\infty} & \sum_{x \in V} u(x, t)\left(\varphi_R\right)_{t t}(x, t) \mu(x) d t \mid \\
	   	 \leqslant & \int_0^{\infty} \sum_{x \in V}|u(x, t)|\left|\left(\varphi_R\right)_{t t}(x, t)\right| \mu(x) d t \\ 
	   	 \leqslant & \frac{C}{R^{1+\alpha}} \int_{R^{\frac{1+\alpha}{2}}}^{2 R^{\frac{1+\alpha}{2}}} \sum_{x \in V}|u(x, t)| \eta^{s-2}\left(\frac{t}{R^{\frac{1+\alpha}{2}}}\right) e^{-\delta \frac{d\left(x, x_0\right)}{R}} \mu(x) d t \\
	   	  \leqslant &\frac{C}{R^{1+\alpha}} \left(\int_{R^{\frac{1+\alpha}{2}}}^{2 R^{\frac{1+\alpha}{2}}} \sum_{x \in V}|u(x, t)|^q h_2(x, t) \eta^s\left(\frac{t}{R^{\frac{1+\alpha}{2}}}\right) e^{-\delta \frac{d\left(x, x_0\right)}{R}} \mu(x) d t\right)^{\frac{1}{q}} \\ 
	   	  & \times\left( \int_{R^{\frac{1+\alpha}{2}}}^{2 R^{\frac{1+\alpha}{2}}} \sum_{x \in V} h_2^{-\frac{1}{q-1}}(x, t) \eta^{s-\frac{2 q}{q-1}}\left(\frac{t}{R^{\frac{1+\alpha}{2}}}\right) e^{-\delta \frac{d\left(x, x_0\right)}{R}} \mu(x) d t\right)^{\frac{q-1}{q}} \\ 
	   	  \leqslant & C\left(\int_{R^{\frac{1+\alpha}{2}}}^{2 R^{\frac{1+\alpha}{2}}} \sum_{x \in V}|u(x, t)|^q h_2(x, t) \mu(x) d t\right)^{\frac{1}{q}}R^{(\alpha+1)\frac{-p+q}{(pq-1)q}}\to 0\mbox{ as }R\to\infty,\end{aligned}
	   \end{equation}
	   where we have used again \eqref{ew7}.
	   Letting $R\to\infty$ in \eqref{eat41} and using \eqref{ew4}, \eqref{ew8},  \eqref{eat51} and \eqref{eat61}, we obtain
	   \begin{equation*}
	   	\int_0^{\infty} \sum_{x \in V} \mu(x) h_1(x,t)|v(x, t)|^p d t\leq 0, 
	   \end{equation*}
	   which implies $v=0$. Similarly, we also obtain $u=0$. The proof is completed.\qed
	   

\bibliographystyle{abbrv}

\end{document}